\documentclass[12pt, reqno]{amsart}
\usepackage{ifthen,amsfonts,amsmath,amssymb,epic,eepic}
\usepackage{epsfig,color}

\makeatother

\theoremstyle{definition}

\newtheorem{Thm}{Theorem}

\newtheorem{Lem}[Thm]{Lemma}

\newtheorem{Def}[Thm]{Definition}

\newtheorem{Exp}{Example}

\def\RR{{\mathbf R}}

\def\MM{{\mathcal{M}}}
\def\AA{{\mathcal{A}}}
\newcommand{\eqr}[1]{(\ref{#1})}

\begin{document}

\title{Bounding dimension of ambient space by density for Mean Curvature Flow}

\author{Maria Calle}
\address{Courant Institute of Mathematical Sciences\\
251 Mercer Street\\ New York, NY 10012}

\email{calle@cims.nyu.edu}

\date{\today}

\maketitle

Abstract: For an ancient solution of the mean curvature flow, we show that each time slice $M_t$ is contained in an affine subspace with dimension bounded in terms of the density and the dimension of the evolving submanifold. Recall that an ancient solution is a family $M_t$ that evolves under mean curvature flow for all negative time $t$.

\section{Introduction} \label{s:s0}
\vskip6mm

This paper deals with ancient solutions of mean curvature flow. An ancient solution is a family $(M_t)$ of $n$-dimensional submanifold of $\RR^{n+k}$ that moves by mean curvature flow for all negative time $t$ (or in general, for all times $t<T$ for some fixed $T$). We prove that each $M_t$ is contained in an affine subspace of bounded dimension. The bound on the dimension depends only on a bound on the density and on the dimension of the evolving manifold.

\vskip6mm

A family $(M_t)_{t\in(a,b)}$ of $n$-dimensional submanifolds of $\RR^{n+k}$ moves by mean curvature if there exist immersions $x_t=x(\cdot,t):M^n\longrightarrow\RR^{n+k}$ of an $n$-dimensional manifold $M^n$ with images $M_t=x_t(M^n)$ satisfying the evolution equation
\begin{equation}  \label{mcf}
\frac{\partial x}{\partial t}=\vec{H}.
\end{equation}
Here $\vec{H}(p,t)$ denotes the mean curvature vector of $M_t$ at $x(p,t)$ for $(p,t)\in M^n\times(a,b)$. The space-time track of the family $(M_t)$ is the set
\begin{equation}
\mathcal{M}=\bigcup_{t\in(a,b)}M_t\times\{t\}\subset\RR^{n+k}\times\RR, \notag
\end{equation} 
sometimes simply written as $\mathcal{M}=\{(M_t,t),t\in(a,b)\}$.
\vskip6mm

In particular, a minimal $n$-dimensional submanifold $M$ of $\RR^{n+k}$ is a stationary solution of the evolution equation $\eqr{mcf}$, because in $M$ we have $\vec{H}=0$.\vskip6mm

Mean curvature flow can be defined not only for smooth manifolds, but also for more general objects. In particular, in \cite{B} Brakke defines mean curvature flow for integral varifolds. A varifold is a measure-theoretic generalization of a manifold that can have singularities. Often a smooth solution of $\eqr{mcf}$ develops singularities in finite time, and after that it becomes a varifold solution (also called a weak solution). Most of what we state in this paper for smooth solutions of mean curvature flow is also true for weak solutions.

\vskip6mm

For a minimal $n$-dimensional submanifold the density at a point $x_0$ is defined by $\Theta(M,x_0)=\lim_{r\to 0^+}\frac{\mathcal{H}^n(M\cap B_r)}{r^n}$. For a solution of \eqr{mcf}, instead of the area $\mathcal{H}^n(M\cap B_r)$ we consider an integral quantity that we denote by $\AA(\MM\cap E_r)$, defined by Ecker in \cite{E1} (see section $\ref{s:s1}$ of this paper for a precise definition). This quantity plays the role of the area $\mathcal{H}^n(M\cap B_r)$, in particular the density for a mean curvature flow solution $\MM$ at a point $(x_0,t_0)\in\MM$ is defined by $\Theta(\MM,x_0,t_0)=\lim_{r\to 0^+}\frac{\AA(\MM\cap E_r)}{r^n}$.

\vskip6mm

Philosophically, one can think of the ratio $\frac{\AA(\MM\cap E_r)}{r^n}$ as a measure for how close the space-time track $\MM$ is to $\RR^{n+1}$. For $\MM=\RR^{n+1}$, this ratio is constantly equal to $1$. As we will see, the main result in this paper goes along this line: we prove that when this ratio is bounded uniformly for all $r$, the manifolds $M_t$ are contained in affine subspaces of dimension possibly smaller than the dimension of the ambient space $\RR^{n+k}$. For a solution of mean curvature flow, this ratio is nondecreasing in $r$ (see section $\ref{s:s1}$), and its limit when $r$ goes to $0$ is called the density. For smooth solutions the density is always $1$. On the other hand, unit density gives regularity: a weak solution with unit density almost everywhere is smooth (see \cite{B}).

\vskip6mm

Our main result is the following:
\begin{Thm}     \label{thm1}
Let $\mathcal{M}$ be an ancient solution of MCF in $\RR^{n+k}$ such that $\forall t\in(-\infty,0)$, $M_t$ has no boundary in $\RR^{n+k}$ and has finite mass $\mathcal{H}^n(M_t\cap B_{2\sqrt{-2nt}})<\infty$. If $\MM$ satisfies:
\begin{equation}
1\le\frac{\mathcal{A}(\mathcal{M}\cap E_r)}{r^n}\le V_{\mathcal{M}}<\infty \quad\forall r>0, \notag
\end{equation}
then for each $t\in(-\infty,0)$, $M_t$ is contained in some affine subspace with
\begin{equation}
\dim\le(n+1)\,\frac{n}{n-1}\,2^{5n+12}\,V_{\mathcal{M}}.\notag
\end{equation}
\end{Thm}
\vskip6mm

The main point of this bound on the dimension is that it only depends on $n$, but not on $k$. We have stated the theorem for smooth solutions of $\eqr{mcf}$, but in fact we will see that it also holds for varifold solutions of mean curvature flow. As mentioned above, this is important since smooth solutions very often develop singularities.
\vskip6mm

In the theorem we ask that the manifolds $M_t$, $t\in(-\infty,0)$ have no boundary in $\RR^{n+k}$ and have finite mass $\mathcal{H}^n(M_t\cap B_{2\sqrt{-2nt}})<\infty$. An ancient solution $\MM$ with this property is called well-defined in $\RR^{n+k}$. Throughout the paper we assume that $\MM$ is a well-defined ancient solution. As seen in \cite{E2}, this guarantees that all integral quantities considered in the paper are finite.
\vskip6mm

The organization of the paper is as follows: in the first section, we recall some facts about mean curvature flow, in particular a mean value formula proved by Ecker in \cite{E2} that will be essential for our proof. In this section, we also develop a little the concept of weak solutions of mean curvature flow, and give some examples of smooth ancient solutions. In section 2, we give the statement of a second theorem, from which our main theorem is a corollary. Also, we give an idea of the steps of the proof. Many of the ideas for this proof are inspired by the paper \cite{CM2} of Colding and Minicozzi. In section 3, we give the statement and proof of several lemmas, necessary for the proof of theorem $\ref{thm2}$. Finally, in section 4 we give the proof of theorem $\ref{thm2}$.
\vskip6mm

I would like to thank my advisor Tobias Colding for his help and for suggesting this problem.
\vskip6mm

\section{Preliminaires}  \label{s:s1}

In this section we introduce some definitions and formulas about mean curvature flow, give an idea of how Brakke solutions of mean curvature flow are defined, and give some examples of ancient solutions.
\vskip6mm

A minimal $n$-dimensional submanifold $M$ on $\RR^{n+k}$ is a stationary solution of the evolution equation $\eqr{mcf}$. For such a manifold $M$, we have the following monotonicity formula (see \cite{CM1}):
\begin{equation}
\frac{d}{dr}\left(\frac{\mathcal{H}^n(M\cap B_r)}{r^n}\right)=\frac{d}{dr}\int_{M\cap B_r}\frac{|x^{\bot}|^2}{|x|^{n+2}}d\mathcal{H}^n.\notag
\end{equation}
In fact this formula is a consequence of the following mean value inequality (see also \cite{CM1}):
\begin{equation}
\frac{d}{dr}\left(\frac{1}{r^n}\int_{B_r\cap M}f\right)=\frac{1}{r^{n+1}}\int_{\partial B_r\cap M}f \frac{|(x-x_0)^N|^2}{|(x-x^0)^{\bot}|}+\frac{1}{2r^{n+1}}\int_{B_r\cap M}(r^2-|x-x_0|^2)\triangle_Mf.\notag
\end{equation}
\vskip6mm

For a solution of the MCF $(M_t)_{t\in(a,b)}$, we have a different monotonocity formula: for $r>0$, we define the `heatball' $E_r=E_r(0,0)$ centered at $(0,0)\in\RR^{n+k}\times\RR$ to be the bounded open set
\begin{equation}
E_r=\{(x,t)\in\RR^{n+k}\times\RR,\; t<0,\; \psi_r(x,t)>0\}\subset\RR^{n+k}\times\RR,\notag
\end{equation}
where $\psi_r\equiv\log{\Phi r^n}\equiv\psi+n\log{r}$ and $\Phi(x,t)=\frac{1}{(-4\pi t)^{\frac{n}{2}}}e^{\frac{|x|^2}{4t}}$. Observe that $E_r$ can also be written as:
\begin{equation}
E_r=\bigcup_{-r^2/4\pi<t< 0}B_{R_r(t)}\times\{t\}, \notag
\end{equation}
where $R_r(t)=\sqrt{2nt\log(-4\pi t/r^2)}$.
Following Ecker (see \cite{E2}), we introduce the integral quantity
\begin{equation} \label{defAA}
\mathcal{A}(\mathcal{M}\cap E_r)\equiv\int\int_{\mathcal{M}\cap E_r}|\nabla\psi|^2+|\vec{H}|^2\psi_r,
\end{equation}
which will play the same role as the area $\mathcal{H}^n(M\cap B_r)$ for minimal submanifolds. Here we use the shorthand notation
\begin{equation}
\int\int_{\mathcal{M}\cap E_r}f\equiv\int^0_{-\frac{r^2}{4\pi}}\int_{M_t\cap B_{R_r(t)}}fd\mu_tdt,\notag
\end{equation}
where $d\mu_t$ denotes the surfaces measure on $M_t$ and $\nabla=\nabla^{M_t}$ denotes the tangential component of the gradient of a function at a point of $M_t$. For notational simplicity, we will set $\Gamma=|\nabla\psi|^2+|\vec{H}|^2\psi_r$.
\vskip6mm

The quantity defined in $\eqr{defAA}$ behaves like $n$-dimensional area in a number of ways. In particular, for a solution $(M_t)_{t\in(a,b)}$ of the MCF Ecker proved the following mean value formula (see \cite{E2}):
\begin{equation} \label{mvf}
\frac{d}{dr}\left(\frac{1}{r^n}\int\int_{\mathcal{M}\cap E_r}u\:\Gamma\right)=\frac{n}{r^{n+1}}\int\int_{\mathcal{M}\cap E_r}\left(u|\vec{H}-\nabla^{\bot}\psi|^2-\psi_r\left(\frac{d}{dt}-\triangle\right)u\right),
\end{equation}
where $u$ is any function for which all integral expressions are finite. Setting $u=1$ in $\eqr{mvf}$, we get that $r^{-n}\AA(\MM\cap E_r)$ is increasing in $r$. Thus the density:
\begin{equation}
\Theta(\MM,0,0)=\lim_{r\to 0}\frac{\AA(\MM\cap E_r)}{r^n} \notag
\end{equation}
at the point $(0,0)$ is well-defined. Moreover, if $u\ge 0$ is continuous at $(0,0)$ and $(\frac{d}{dt}-\triangle)u\le 0$, then $\eqr{mvf}$ implies:
\begin{equation} \label{mvi}
u(0,0)\,\Theta(\mathcal{M},0,0)\le\frac{1}{r^n}\int\int_{\mathcal{M}\cap E_r}u\,\Gamma.
\end{equation}
\vskip6mm

By replacing $\psi$ by $\psi_{x_0,t_0}(x,t)\equiv\psi(x-x_0,t-t_0)$ we obtain analogous statements for MCF in heat-balls $E_r(x_0,t_0)$, the translates of $E_r$. In particular, we get that the quantity $r^{-n}\AA_{x_0,t_0}(\MM\cap E_r(x_0,t_0))$ is nondecreasing in $r$, where $\AA_{x_0,t_0}(\MM\cap E_r(x_0,t_0))$ is $\AA(\MM\cap E_r)$ translated to $(x_0,t_0)$. For notational convenience, all results are stated with respect to the reference point $(0,0)\in\RR^{n+k}\times\RR$.
\vskip6mm

So far we have been talking about smooth solutions of mean curvature flow. In order to introduce the concept of weak solutions of the differential equation $\eqr{mcf}$, we consider an integral version of this equation, which first appeared in \cite{B}. This integral version is stated by Ecker in \cite{E1} as follows: for any smooth solution $(M_t)_{t\in (a,b)}$ of mean curvature flow in an open subset $U\in\RR^{n+1}$ we have:
\begin{equation} \label{intmcf}
\frac{d}{dt}\int_{M_t}\phi=\int_{M_t}\vec{H}\cdot D\phi-|\vec{H}|^2\phi
\end{equation}
for all $t\in(a,b)$ and $\phi\in C^1_0(U)$. This integral equation can be used as a definition of mean curvature flow since any family of smooth submanifolds satisfying $\eqr{intmcf}$ is also a solution of $\eqr{mcf}$. 
\vskip6mm

The equation $\eqr{intmcf}$ is the motivation of Brakke solutions (see \cite{B}). Brakke defines mean curvature flow for generalized submanifolds (so-called integral varifolds) in the language of geometric measure theory. In short, an $n$-dimensional varifold in $\RR^{n+k}$ is a Radon measure on $G_n(\RR^{n+k})=\RR^{n+k}\times G(n+k,n)$, where $G(n+k,n)$ is the Grassmann manifold of $n$-dimensional planes of $\RR^{n+k}$. The definition of integral varifolds is more restrictive, so that it allows to define integrals of functions over varifolds. That way an equation like $\eqr{intmcf}$ makes sense for varifolds. One has to allow for sudden local loss of area in this setting so the integral identity $\eqr{intmcf}$ has to be replaced by an inequality of the form LHS$\le$RHS. Moreover, the right-hand side of $\eqr{intmcf}$ requires re-interpretation for integral varifolds, i.e., one has to define the mean curvature vector $\vec{H}$ on a varifold. We won't give here details of this construction, which was first done by Brakke in \cite{B}. For our purposes, it is enough to remark that the mean value formula $\eqr{mvf}$ that we use to prove our result can be derived by substituting appropriate test functions into identity $\eqr{intmcf}$. In this case, the identity $\eqr{mvf}$ has to be replaced by an inequality of the form LHS$\ge$RHS. For the purpose of this paper, this inequality is enough to prove our theorem.
\vskip6mm

Finally, we close this section with some examples of ancient solutions of mean curvature flow:
\vskip6mm

\begin{Exp} 
{\it Minimal submanifolds}. As mentioned above, if $M$ is a minimal submanifold of $\RR^{n+k}$, then $M_t=M$ for all $t\in\RR$ is a solution of the evolution equation $\eqr{mcf}$ defined for all times $t\in(-\infty,\infty)$.
\end{Exp}
\vskip6mm

\begin{Exp}
{\it Shrinking spherical cylinders}. Shrinking spherical cylinders are defined by:
\begin{equation}
M_t=\partial B^{n+1-l}_{r(t)}\times\RR^l \notag
\end{equation}
for $0\le l\le n$ (the case $l=0$ corresponds to shrinking spheres). For the $n$-dimensional submanifolds $M_t\subset\RR^{n+1}$ to satisfy $\eqr{mcf}$, the radius $r(t)$ has to be a solution of the ODE:
\begin{equation}
r'(t)=-\frac{n-l}{r}. \notag
\end{equation}
If we fix the radius at time $0$ to be $r(0)=\rho$, we obtain:
\begin{equation}
r(t)=\sqrt{\rho^2-2(n-l)t} \notag
\end{equation}
and the solution exists for $t\in(-\infty,\frac{\rho^2}{2(n-l)})$. For time $t=\frac{\rho^2}{2(n-l)}$ it degenerates to the $l$-dimensional plane $\{ 0 \}\times\RR^l\subset\RR^{n+1}$, and after that time it disappears.
\end{Exp}
\vskip6mm

\begin{Exp}
{\it Grim reaper}. The grim reaper is a graph solution of MCF, i.e. it has the form $M_t=\mbox{graph}\, u(.,t)$ for $u(.,t)\, :\RR\to\RR$ given by:
\begin{equation}
u(p,t)=-\log\cos{p}+t,\quad p\in(-\frac{\pi}{2},\frac{\pi}{2}). \notag
\end{equation}
In this case, $M_t$ is defined and smooth for all times $t\in(-\infty,\infty)$, and moves by translation.
\end{Exp}
\vskip6mm

Observe that the theorem is not useful for these examples, because here the submanifolds $M_t$ are hypersurfaces of $\RR^{n+1}$.
\vskip6mm

\section{Structure of the proof}   \label{s:s2}
\vskip6mm

In this section we state the theorem that will imply our main result, and give a sketch of the proof.
\vskip6mm

Our main result is a consequence of a theorem that bounds the dimension of certain spaces of functions. Namely, we define the space $\mathcal{H}_d(\mathcal{M})$ to be the linear space of functions $u:\RR^{n+k}\times\RR\longrightarrow\RR$ satisfying
\begin{equation} \label{eq2}
\left(\frac{d}{dt}-\triangle\right)u=0
\end{equation}
and such that $|u(x,t)|\le C(1+|x|^d)$ for all $(x,t)\in\MM$ and for some constant $C<\infty$ (i.e., $u$ has polynomial growth of order at most $d$). With this definition, the coordinate functions $x_i$ are in $\mathcal{H}_1(\mathcal{M})$, because $\frac{d}{dt}x_i=\frac{\partial}{\partial t}x_i+Dx_i\cdot\vec{H}=\vec{H}_i=\triangle_{M_t}x_i$, where $D$ stands for differentiation in $\RR^{n+k}$.
The result that gives us the main theorem is the following
\begin{Thm} \label{thm2}
If $\mathcal{M}$ is well-defined in $\RR^{n+k}\times(-\infty,0)$ satisfying 
\begin{equation}
1\le\frac{\mathcal{A}(\mathcal{M}\cap E_r)}{r^n}\le V_{\mathcal{M}}<\infty \quad\forall r>0, \notag
\end{equation}
then for any $d\ge 1$,
\begin{equation}
\dim\mathcal{H}_d(\mathcal{M})\le C V_{\mathcal{M}}d^{n-1} \notag
\end{equation}
where $C=(n+1)\frac{n}{n-1}2^{5n+12}$.
\end{Thm}
\vskip6mm

To prove this bound on the dimension, we define a family of inner products for functions $u,v$ in L$^2$($\RR^{n+k}\times\RR$) by:
\begin{equation}
J_r(u,v)=\int\int_{\mathcal{M}\cap E_{\rho(r)}}u\,v\,\Gamma.
\end{equation}
Here $r>0$ and $\rho$ is given in terms of $r$ by $\rho(r)=\sqrt{\frac{2\pi e}{n}}r$. We also define, for $u\in$L$^2$($\RR^{n+k}\times\RR$), the following function of $r$:
\begin{equation}
I_u(r)=\int\int_{\mathcal{M}\cap E_{\rho(r)}}u^2\,\Gamma.
\end{equation}
Observe that the function $I_u(r)$ correspond to the quadratic form associated to $J_r$, i.e., $I_u(r)=J_r(u,u)$. If $u$ has polynomial growth of order at most $d$, and in particular if $u\in\mathcal{H}_d(\mathcal{M})$, then $I_u(r)\le C(1+r^{2d+n}),\,\forall r>0$.
\vskip6mm

The proof of theorem $\ref{thm2}$ has then two steps: first, we choose $\Omega=(1-\frac{1}{2d})^{-1}$. Then, given linearly independent functions on $\mathcal{H}_d(\MM)$, we find a family of $J_r$-orthonormal functions (where $r=\Omega^{m+1}$ for a certain integer $m$) satisfying $(\frac{d}{dt}-\triangle)v=0$ with a uniform lower bound on $I_v(\Omega^{-1}r)$.
\vskip6mm

The second step is to reduce the problem to bounding the number of $J_r$-orthonormal functions as above. To do so, we construct a function $K(x,t)=\sum_{i=1}^L|v_i(x,t)|^2$, where $v_1,\ldots,v_L$ are orthonormal functions. We give a bound of $I_K(\Omega^{-1}r)$ which depends on $\Omega$ but not on $L$, and using the lower bound of $I_{v_i}(\Omega^{-1}r)$ we get a bound on $L$ that depends on $\Omega$. This bound gives a bound in the number of linearly independent functions with which we started, hence giving a bound on the dimension of $\mathcal{H}_d(\MM)$. The bound turns out to be polynomial in $d$.
\vskip6mm

As mentioned above, these ideas (the definition of the functions $J_r, \, I_u$ and of the ``Bergman kernel'' $K$, the steps of the proof) are inspired by the main theorem proved by Colding and Minicozzi in \cite{CM2}.
\vskip6mm

\section{General constructions}  \label{s:s3}

In this section we state and proof some technical lemmas necessary for the proof.
\vskip6mm

The following lemma is stated and proven in \cite{CM1}. To make things self-contained, we include the proof.
\begin{Lem} \label{lemma1}
Suppose that $f_1,\ldots,f_{2s}$ are nonnegative nondecreasing functions on $(0,\infty)$ such that none of the $f_i$ vanishes identically, and for some $d_0,K>0$ and all $i$:
\begin{equation}
f_i(r)\le K(r^{d_0}+1)\quad   \forall r>0. \notag
\end{equation}
Then for all $\Omega >1$, there exist $s$ of these functions, $f_{\alpha_1},\ldots,f_{\alpha_s}$, and infinitely many integers, $m\ge 1$, such that for $i=1,\ldots,s$
\begin{equation}
f_{\alpha_i}(\Omega^{m+1})\le\Omega^{2d_0}f_{\alpha_i}(\Omega^m).  \notag
\end{equation}
\end{Lem}
\vskip6mm

\begin{proof}  \label{pflemma1}
Since the functions are nondecreasing and none of them vanishes identically, we may suppose that for some $R>0$ and any $r\ge R,\quad f_i(r)>0$ for all $i$. We will show that there are infinitely many $m$ such that there is some rank $s$ subset of $\{f_i\}$ (where the subset could vary with $m$) satisfying the inequality. This will suffice to prove the lemma, since there are only finitely many rank $s$ subsets of the set $\{f_1,\ldots,f_{2s}\}$, hence one of these rank $s$ subsets must be repeated infinitely often.
For $r>R$, note that:
\begin{equation}  \label{pfl1eq1}
g(r)=\prod_{i=1}^{2s}f_i(r)\le K^{2s}(r^{d_0}+1)^{2s}
\end{equation}
and $g$ is a positive nondecreasing function. Assume that there are only finitely may $m\ge\frac{\log{R}}{\log{\Omega}}$ satisfying the inequality for at least $s$ of the functions $f_i$. Let $m_0-1$ be the largest such $m$. For all $j\ge 1$ we have that
\begin{equation}
\Omega^{2d_0(s+1)}g(\Omega^{m_0+j-1})<g(\Omega^{m_0+j}).  \notag
\end{equation}
Iterating this and applying $\eqr{pfl1eq1}$ gives for any $j\ge 1$
\begin{equation}
\Omega^{2d_0(s+1)j}g(\Omega^{m_0})<g(\Omega^{m_0+j})\le C(\Omega^j)^{2sd_0},  \notag
\end{equation}
where $C=C(s,m_0,\Omega,K,d_0)$. Since $\Omega>1$, taking $j$ large yields the contradiction.
\end{proof}
\vskip6mm

Given a linearly independent set of functions in $\mathcal{H}_d(\mathcal{M})$, we will construct functions of one variable which reflect the growth and independence properties of this set. We begin with two definitions:
\vskip6mm

\begin{Def} \label{def1}
Suppose that $u_1,\ldots,u_s$ are linearly independent functions on $\mathcal{H}_d(\MM)$. For each $r>0$ we will now define a $J_r$-orthogonal spanning set $w_{i,r}$ and functions $f_i$, $i=1,\ldots,s$. Set $w_{1,r}=u_1$ and $f_1(r)=I_{w_{1,r}}(r)$. Define $w_{i,r}$ by requiring it to be orthogonal to $u_j$ for $j<i$ with respect to the inner product $J_r$ and so that
\begin{equation}
u_i=\sum_{j=1}^{i-1}\lambda_{ji}(r)u_j+w_{i,r}.
\end{equation}
Note that $\lambda_{ji}(r)$ is not uniquely defined if the $u_i|E_{\rho(r)}$ are linearly dependent. However, since the $u_i$ are linearly independent, $\lambda_{ji}(r)$ will be uniquely defined for $r$ sufficiently large. In any case, the following quantity is well-defined for all $r>0$:
\begin{equation}
f_i(r)=\int\int_{\mathcal{M}\cap E_{\rho(r)}}w_{i,r}^2\,\Gamma=I_{w_{i,r}}(r).  \notag
\end{equation}
\end{Def}
\vskip6mm

In the next lemma, we will record some properties of the functions $f_i$:
\begin{Lem} \label{lemma2}
If $u_1,\ldots,u_s\in\mathcal{H}_d(\mathcal{M})$ are linearly independent,then there exists a constant $K>0$ (depending on the set $\{u_i\}$) such that for $i=1,\ldots,s$, the functions $f_i$ defined above verify:
\begin{equation}   \label{l2eq1}
f_i(r)\le K(r^{2d+n}+1),
\end{equation}
\begin{equation}   \label{l2eq2}
f_i\quad \mbox{is a nondecreasing function,}
\end{equation}
\begin{equation}   \label{l2eq3}
f_i\quad \mbox{is nonnegative, and positive for $r$ sufficently large, and}
\end{equation}
\begin{equation}   \label{l2eq4}
f_i(r)=I_{w_{i,r}}(r) \quad\mbox{and}\quad f_i(t)\le I_{w_{i,r}}(t) \quad\mbox{for}\quad t<r.
\end{equation}
\end{Lem}
\vskip6mm

\begin{proof}  \label{pflemma2}
For $\eqr{l2eq1}$, we observe that:
\begin{eqnarray}
f_i(r) & = & \int\int_{\MM\cap E_{\rho(r)}}w_{i,r}^2\,\Gamma \notag \\
& \le & \int\int_{\MM\cap E_{\rho(r)}}w_{i,r}^2\,\Gamma+\int\int_{\MM\cap E_{\rho(r)}}\left(\sum_{j<i}\lambda_{j,i}(r)u_j\right)^2\,\Gamma \\
& = & \int\int_{\MM\cap E_{\rho(r)}}\left(\sum_{j<i}\lambda_{j,i}(r)u_j+w_{i,r}\right)^2\,\Gamma \notag \\
& = & I_{u_i}(r)\le K(r^{2d+n}+1), \notag
\end{eqnarray}
where the second equality follows from the orthogonality of $w_{i,r}$ and $u_j$ for $j<i$.
Furthermore, for $s<r$:
\begin{eqnarray} \label{pfl2eq1}
f_i(s) & = & \int\int_{\MM\cap E_{\rho(s)}}|u_i-\sum_{j=1}^{i-1}\lambda_{ji}(s)u_j|^2\,\Gamma=I_{w_{i,s}}(s) \notag \\
& \le & \int\int_{\MM\cap E_{\rho(s)}}|u_i-\sum_{j=1}^{i-1}\lambda_{ji}(r)u_j|^2\,\Gamma=I_{w_{i,r}}(s) \\
& \le & \int\int_{\MM\cap E_{\rho(r)}}|u_i-\sum_{j=1}^{i-1}\lambda_{ji}(r)u_j|^2\,\Gamma=I_{w_{i,r}}(r)=f_i(r). \notag
\end{eqnarray}
Here the first inequality follows from:
\begin{eqnarray}
I_{w_{i,r}}(s)-I_{w_{i,s}}(s) & = & \int\int_{\MM\cap E_{\rho(s)}}(w_{i,r}+w_{i,s})(w_{i,r}-w_{i,s})\,\Gamma \notag \\
& = & \int\int_{\MM\cap E_{\rho(s)}}(w_{i,r}+w_{i,s})\left(\sum_{j<i}(\lambda_{ji}(s)-\lambda_{ji}(r))u_j\right)\,\Gamma \\
& = & \int\int_{\MM\cap E_{\rho(s)}}w_{i,r}(w_{i,r}-w_{i,s})\,\Gamma \notag \\
& = & \int\int_{\MM\cap E_{\rho(s)}}(w_{i,r}-w_{i,s})^2\,\Gamma\ge 0, \notag
\end{eqnarray}
and the second inequality in $\eqr{pfl2eq1}$ follows from the monotonicity of $I_{w_{i,r}}$. The inequalities in $\eqr{pfl2eq1}$ imply $\eqr{l2eq4}$, and from them and the linear independence of the $u_i$, we get $\eqr{l2eq2}$ and $\eqr{l2eq3}$.
\end{proof}
\vskip6mm

In the next lemma, we apply lemma \ref{lemma1} to the functions $f_i$:
\vskip6mm

\begin{Lem} \label{lemma3}
Suppose that $u_1,\ldots,u_{2s}\in\mathcal{H}_d(\mathcal{M})$ are linearly independent. Given $\Omega >1$ and $m_0>0$, there exist $m\ge m_0$, an integer $L\ge\frac{1}{2}\Omega^{-4d-2n}s$, and functions $v_1,\ldots,v_L$ in the linear span of the $u_i$ such that for $i,j=1,\ldots,L$
\begin{equation}
J_{\Omega^{m+1}}(v_i,v_j)=\delta_{i,j}  \notag
\end{equation}
and
\begin{equation}
\frac{1}{2}\Omega^{-4d-2n}\le I_{v_i}(\Omega^m).  \notag
\end{equation}
\end{Lem}
\vskip6mm

\begin{proof}  \label{pflemma3}
By lemma $\ref{lemma2}$, we can apply lemma $\ref{lemma1}$ to the functions $f_i$ in definition $\ref{def1}$ with $d_0=2d+n$. Therefore there exist $m\ge m_0$ and a subset $f_{\alpha_1},\ldots,f_{\alpha_s}$ such that for $i=1,\ldots,s$:
\begin{equation}  \label{pfl3eq1}
0<f_{\alpha_i}(\Omega^{m+1})\le\Omega^{4d+2n}f_{\alpha_i}(\Omega^m).
\end{equation}
Let $w_{\alpha_i,\Omega^{m+1}},\:i=1,\ldots,s$ be the corresponding functions in the linear span of the $u_i$ as in definition $\ref{def1}$. Consider the $s$-dimensional linear space spanned by the functions $w_{\alpha_i,\Omega^{m+1}}$ with inner product $J_{\Omega^{m+1}}$. On this space there is also the positive semidefinite bilinear form $J_{\Omega^m}$. Let $v_1,\ldots,v_s$ be an orthonormal basis for $J_{\Omega^{m+1}}$ which diagonalizes $J_{\Omega^m}$. We will now evaluate the trace of $J_{\Omega^m}$ with respect to these two bases. First, with respect to the orthogonal basis $w_{\alpha_i,\Omega^{m+1}}$ we get, by $\eqr{pfl3eq1}$ and $\eqr{l2eq4}$:
\begin{equation}
s\,\Omega^{-4d-2n}\le\sum_{i=1}^s\frac{I_{w_{\alpha_i,\Omega^{m+1}}}(\Omega^m)}{I_{w_{\alpha_i,\Omega^{m+1}}}(\Omega^{m+1})}.  \notag
\end{equation}
Since the trace is independent of the choice of basis we get when evaluating this on the orthonormal basis $v_i$:
\begin{equation}
s\,\Omega^{-4d-2n}\le\sum_{i=1}^sI_{v_i}(\Omega^m).  \notag
\end{equation}
Combining this with $0\le I_{v_i}(\Omega^m)\le I_{v_i}(\Omega^{m+1})=1$, which follows from the monotonicity of $I$, we get that there exist at least $l\ge\frac{s}{2}\Omega^{-4d-2n}$ of the $v_i$ such that:
\begin{equation}
\frac{1}{2}\Omega^{-4d-2n}\le I_{v_i}(\Omega^m)\le I_{v_i}(\Omega^{m+1})=1.  \notag
\end{equation}
This proves the lemma.
\end{proof}
\vskip6mm

We have defined the quantity $\mathcal{A}(\mathcal{M}\cap E_r)$ only for the heat-ball $E_r$, but it can be defined similarly for a general set $C\subset\RR^{n+k}$ by setting $\mathcal{A}(\mathcal{M}\cap C)=\int\int_{\mathcal{M}\cap C}\Gamma$. We have then the following lemma:
\begin{Lem} \label{lemma4}
For $\mathcal{M}$ as in theorem $\ref{thm2}$, for all $0<\delta\le\frac{1}{2}$, there is a $R_0=R_0(\delta)>0$ with
\begin{equation}
\mathcal{A}(\mathcal{M}\cap E_r\backslash E_{(1-\delta)r})\le 2(n+1)\,\delta\,\mathcal{A}(\mathcal{M}\cap E_r) \quad\forall r\ge R_0.  \notag
\end{equation}
\end{Lem}
\vskip6mm

\begin{proof}  \label{pflemma4}
The mean value formula gives that $\frac{\AA(\MM\cap E_r}{r^n}$ is nondecreasing in $r$, and it is bounded for all $r$ by hypothesis, hence we can define $V_{\MM}$ by:
\begin{equation}
V_{\MM}=\lim_{r\to\infty}\frac{\AA(\MM\cap E_r)}{r^n}<\infty. \notag
\end{equation}
This implies that, given $0<\delta\le\frac{1}{2}$, we can choose $R_0$ so that $\forall R\ge\frac{R_0}{2}$,
\begin{equation} \label{pfl4eq1}
V_{\MM}-\frac{\AA(\MM\cap E_R)}{R^n}<\delta V_{\MM}. 
\end{equation}
Now, $\forall R\ge R_0$, apply $\eqr{pfl4eq1}$ to $(1-\delta)R\ge R_0/2$ to get:
\begin{equation}
\AA(\MM\cap E_{(1-\delta)R})\ge(1-\delta)^{n+1}R^nV_{\MM}.  \notag
\end{equation}
This implies
\begin{eqnarray}
\AA(\MM\cap E_R\backslash E_{(1-\delta)R}) & = & \int\int_{\MM\cap E_R}\Gamma - \int\int_{\MM\cap E_{(1-\delta)R}}\Gamma \notag \\
& \le & V_{\MM}R^n-(1-\delta)^{n+1}R^nV_{\MM}  \notag \\ 
& \le & 2(n+1)\,\delta\,\AA(\MM\cap E_R). \notag
\end{eqnarray}
\end{proof}
\vskip6mm

\section{Proof of finite dimensionality and proof of theorem 2}  \label{s:s4}

In this section we prove theorem $\ref{thm2}$, which gives theorem $\ref{thm1}$ as a corollary.
\vskip6mm

Recall that an ancient solution $\MM$ of mean curvature flow is said to be well-defined in $\RR^{n+k}\times(-\infty,0)$ if the manifolds $M_t,\, t\in(-\infty,0)$, have no boundary in $\RR^{n+k}$ and have finite mass $\mathcal{H}^n(M_t\cap B_{2\sqrt{-2nt}})<\infty$.
\vskip6mm

The following lemma combined with lemmas $\ref{lemma3}$ and $\ref{lemma4}$ will give the proof of theorem $\ref{thm2}$:
\begin{Lem} \label{lemma5}
Let $\mathcal{M}$ be well-defined in $\RR^{n+k}\times(-\infty,0)$ with density at least $1$ and
\begin{equation}
\lim_{r\to\infty}\frac{\AA(\MM\cap E_r)}{r^n}=V_{\MM}<\infty.  \notag
\end{equation}
Suppose $d\ge 1,\: R_0$ satisfy:
\begin{equation} \label{l5eq1}
\AA(\MM\cap E_R\backslash E_{(1-\delta)R})\le 2(n+1)\,\delta\,\AA(\MM\cap E_R)  \quad \forall\,\delta\ge\frac{1}{4d},\quad  \forall R\ge R_0.
\end{equation}
Let $0<a<1,\: r>2R_0$ be fixed, and let $v_1,\ldots,v_L$ be $J_r$-orthonormal functions satisfying $(\frac{d}{dt}-\triangle)v_i=0$ and so that:
\begin{equation} \label{l5eq2}
a\le I_{v_i}((1-\frac{1}{2d})r)\quad \forall i=1,\ldots,L.
\end{equation}
Then:
\begin{equation}
L\le Cd^{n-1},  \notag
\end{equation}
where $C=V_{\MM}2^{n+1}(n+1)\frac{n}{n-1}a^{-1}$.
\end{Lem}
\vskip6mm

\begin{proof}  \label{pflemma5}
Because $d\ge 1$, we can choose an integer $N$ such that $d\le N\le 2d$.
Define the following function on $E_r$:
\begin{equation}
K(x,t)=\sum_{i=1}^L|v_i(x,t)|^2.  \notag
\end{equation}
Observe that $K$ is the trace of the symmetric bilinear form $(u,v)\longmapsto <u,v>(x,t)$ for any $u,v$ in the span of the $v_i$. We can always diagonalize such a bilinear form, therefore, given $(x,t)\in E_r$, we can choose a new $J_r$-orthonormal basis $\{w_i\}$ of span($v_i$) such that $w_i(x,t)=0,\quad \forall i=2,\ldots,L$. The trace of a matrix is invariant under orthogonal change of basis, hence we have:
\begin{equation} \label{pfl5eq1}
K(x,t)=\sum_{i=1}^L|w_i(x,t)|^2=|w_1(x,t)|^2.
\end{equation}
Now the functions $w_i\in\mbox{span}\{v_i\}$, hence they also satisfy $(\frac{d}{dt}-\triangle)w_i=0$. This implies
\begin{equation}
\left(\frac{d}{dt}-\triangle\right) w_i^2=-2|\nabla w_i|^2\le0.  \notag
\end{equation}
Therefore we can apply the mean value inequality $\eqr{mvi}$ to the functions $w_i^2$. For $0<s\le 1$ and $(x,t)\in E_{(1-s)r}(0,0)$, it holds:
\begin{eqnarray} \label{pfl5eq2}
|w_i(x,t)|^2 & \le & \frac{1}{s^nr^n}\int\int_{\MM\cap E_{sr}(x,t)}w_i^2\,\Gamma \\
& \le & \frac{1}{s^n}\frac{V_{\MM}}{\AA(\MM\cap E_r(0,0))}, \notag
\end{eqnarray}
where in the second inequality we have used that $E_{sr}(x,t)\subset E_r(0,0),\quad J_r(w_i,w_i)=1$ and $\frac{\AA(\MM\cap E_r(0,0))}{r^n}\le V_{\MM}$.
Combining $\eqr{pfl5eq1}$ and $\eqr{pfl5eq2}$ we get, for all $(x,t)\in E_{(1-\frac{j}{2d})r}, j=1,\ldots,N$:
\begin{equation} \label{pfl5eq3}
K(x,t)=|w_1(x,t)|^2\le V_{\MM}\left(\frac{j}{2d}\right)^{-n}\frac{1}{\AA(\MM\cap E_r)}\, .
\end{equation}
We want to bound the integral of $K\,\Gamma$ over $E_{(1-\frac{1}{2d})r}$, so we break it down in two terms:
\begin{equation} \label{pfl5eq4}
\int\int_{\MM\cap E_{(1-\frac{1}{2d})r}}K\,\Gamma = \int\int_{\MM\cap E_{(1-\frac{N}{2d})r}}K\,\Gamma + \sum_{j=1}^{N-1}\int\int_{\MM\cap (E_{(1-\frac{j}{2d})r}\backslash E_{(1-\frac{j+1}{2d})r)}}K\,\Gamma.
\end{equation}
For the first term of $\eqr{pfl5eq4}$, we observe that $1-\frac{N}{2d}\le\frac{1}{2}$, therefore we have:
\begin{equation}
\int\int_{\MM\cap E_{(1-\frac{N}{2d})r}}K\,\Gamma\le V_{\MM}\left(\frac{1}{2}\right)^{-n}\frac{\AA(\MM\cap E_{\frac{r}{2}})}{\AA(\MM\cap E_r)} \le V_{\MM}2^n.
\end{equation}
For the second term, we have:
\begin{eqnarray} \label{pfl5eq6}
\sum_{j=1}^{N-1}\int\int_{\MM\cap (E_{(1-\frac{j}{2d})r}\backslash E_{(1-\frac{j+1}{2d})r})}K\,\Gamma & \le & V_{\MM}\sum_{j=1}^{N-1}\left(\frac{j}{2d}\right)^{-n}\frac{\AA(\MM\cap(E_{(1-\frac{j}{2d})r}\backslash E_{(1-\frac{j+1}{2d})r}))}{\AA(\MM\cap E_r)} \notag \\
& \le & V_{\MM}\,2(n+1)\left(\frac{1}{2d}\right)^{-n+1}\,\sum_{j=1}^{N-1}j^{-n} \\
& \le & V_{\MM}\,(n+1)\,\frac{n}{n-1}\,2^nd^{n-1}, \notag
\end{eqnarray}
where the second inequality in $\eqr{pfl5eq6}$ follows from $\eqr{l5eq1}$ and the fact that $E_{(1-\frac{j}{2d})r}\subset E_r$ for all $j=1,\ldots,N-1$, and the third inequality in $\eqr{pfl5eq6}$ follows from the elementary fact:
\begin{equation}
\sum_{j=1}^{N-1}j^{-n}=1+\sum_{j=2}^{N-1}j^{-n}\le 1+\int_1^{\infty}s^{-n}ds=\frac{n}{n-1}\, .  \notag
\end{equation}
This gives the bound:
\begin{eqnarray} \label{pfl5eq5}
\int\int_{\MM\cap E_{(1-\frac{1}{2d})r}}K\,\Gamma & \le & V_{\MM}\,2^n(1+(n+1)\frac{n}{n-1}d^{n-1}) \\
& \le & V_{\MM}\,2^{n+1}(n+1)\frac{n}{n-1}d^{n-1}. \notag
\end{eqnarray}
Finally, by $\eqr{pfl5eq5}$ and $\eqr{l5eq2}$ we get:
\begin{eqnarray}
L\,a & \le & \sum_{i=1}^L I_{v_i}((1-\frac{1}{2d})r) = \int\int_{\MM\cap E_{(1-\frac{1}{2d})r}}K\,\Gamma  \notag \\
& \le & V_{\MM}\,2^{n+1}(n+1)\,\frac{n}{n-1}\,d^{n-1}. \notag
\end{eqnarray}
\end{proof}
\vskip6mm

We are now ready to prove theorem $\ref{thm2}$:

\begin{proof}[Proof of Theorem 2]  \label{Thm2}
Suppose $u_1,\ldots,u_{2s}\in\mathcal{H}_d(\MM)$ are linearly independent. Set $\Omega=(1-\frac{1}{2d})^{-1}>1$, and choose $m_0$ such that $\Omega^{m_0}\ge R_0$, where $R_0=R_0(\delta)$ is chosen as in lemma $\ref{lemma4}$ for $\delta=1/4d$. Then, if we set  $r=\Omega^{N+1}$, we have by lemma $\ref{lemma3}$ that there exists an integer $N\ge m_0$, and $J_r$-orthonormal functions $w_1,\ldots,w_L\in\mbox{span}(u_i)$ such that:
\begin{equation}
\frac{1}{2^{2n+5}}\le I_{w_i}((1-\frac{1}{2d})r)  \notag
\end{equation}
for all $i=1,\ldots,L$. Here $L$ is an integer with
\begin{equation}
L\ge\Omega^{-4d-2n}\frac{s}{2}\ge 2^{-2n-5}s,  \notag
\end{equation}
also by lemma $\ref{lemma3}$. Now we can apply lemma $\ref{lemma5}$ with $a=2^{-2n-5}$, and we get:
\begin{equation}
L\le V_{\MM}(n+1)\,\frac{n}{n-1}\,2^{3n+6}d^{n-1}.  \notag
\end{equation}
Now, because $2s\le 2^{2n+6}L$, we have the desired bound on the dimension:
\begin{equation}
\dim{\mathcal{H}_d(\MM)}\le V_{\MM}(n+1)\frac{n}{n-1}\,2^{5n+12}d^{n-1}.  \notag
\end{equation}
\end{proof}


\begin{thebibliography}{999}
\frenchspacing
\bibitem[B]{B}
K. A. Brakke, The Motion of a Surface by its Mean Curvature, Math. Notes Princeton, NJ, Princeton University Press, 1978
\bibitem[CM1]{CM1}
T.H. Colding and W.P. Minicozzi II, Minimal surfaces, Courant
Lecture Notes in Math., v. 4, 1999.
\bibitem[CM2]{CM2}
T.H. Colding and W.P. Minicozzi II, Liouville theorems for harmonic sections and applications, {\it Communications on Pure and Appied Mathematics} 51 (1998), n. 2, 113-138
\bibitem[E1]{E1}
K. Ecker, Regularity Theory for Mean Curvature Flow, Birkhäuser, 2004.
\bibitem[E2]{E2}
K. Ecker, A local monotonicity formula for mean curvature flow, {\it Annals of Mathematics} 154 (2001), 503-525

\end{thebibliography}
\end{document}